\def\titlep{Classification of sectors of the Cuntz algebras
by graph invariants}
\documentclass[11pt]{article}
\usepackage{graphicx,ifthen}

\font\germ=eufm10 at12pt

\def\goth#1{\hbox{\germ#1}}

\newcommand{\qed}{\hbox{\rule[-2pt]{3pt}{6pt}}}
\newcommand{\qedh}{\hfill\qed \\}
\newcommand{\vv}{\vspace{.3in}}

\newtheorem{Thm}{Theorem}[section]

\newtheorem{defi}[Thm]{Definition}

\newtheorem{prop}[Thm]{Proposition}

\newcommand{\ww}{\vv\noindent}

\def\cal#1{\mathcal #1}
\def\con{{\cal O}_{N}}
\def\edot{=1,\ldots,N}
\def\pr{{\it Proof.}\quad}

\def\co#1{{\cal O}_{#1}}
\addtocounter{footnote}{1}
\def\sftt#1{
\setcounter{equation}{0}
\addtocounter{footnote}{1}
\section{#1}
}

\def\ssft#1{\subsection{#1}}

\def\cls{\quad
\clearpage
}

\begin{document}
%
%
\def\emailp{e-mail:kawamura@kurims.kyoto-u.ac.jp.}
\title{{\bf \titlep}}
\author{\autherp\footnote{\emailp}\\
{\it \addressp}}
\def\inputcls#1{\cls\input #1.txt}
\def\plan#1#2{\par\noindent\makebox[.5in][c]{#1}
\makebox[.1in][l]{$|$}
\makebox[3in][l]{#2}\\}
\def\nset#1{\{1,\ldots,N\}^{#1}}
\def\bfsnl{{\rm BFS}_{N}(\Lambda)}
\def\blank{\cls \quad}
\def\pline{
\noindent
\thicklines
\begin{picture}(1000,5)
\put(0,0){\line(1,0){4250}}
\end{picture}
\par
}
\newcommand{\mline}{\noindent
\thicklines
\begin{picture}(1000,5)
\put(0,0){\line(1,0){750}}
\end{picture}
}
\newcommand{\nline}{\noindent
\thicklines
\begin{picture}(1000,5)
\put(0,0){\line(1,0){3060}}
\end{picture}
\\
}
\def\enda{{\rm End}{\cal A}}
\def\endcon{{\rm End}\con}

\def\autherp{Katsunori  Kawamura}
\def\emailp{e-mail:kawamura@kurims.kyoto-u.ac.jp.}
\def\addressp{College of Science and Engineering Ritsumeikan University,\\
1-1-1 Noji Higashi, Kusatsu, Shiga 525-8577,Japan
}

\def\brl{branching law}

%
\pagestyle{plain}
\setcounter{page}{1}
\setcounter{section}{0}
\setcounter{footnote}{0}

%
%
\begin{center}
{\Large \titlep}

\ww
\autherp
\footnote{\emailp}

\noindent
{\it \addressp}
\quad \\

\end{center}


%
%
\begin{abstract}
A unitary equivalence class of endomorphisms 
of a unital C$^{*}$-algebra ${\cal A}$ is called a {\it sector} of ${\cal A}$.
We introduced permutative endomorphisms
of the Cuntz algebra $\con$ in the previous work.
Branching laws of permutative representations of $\con$
by them are computed by directed regular graphs.
In this article, we classify sectors associated with permutative endomorphisms 
of $\con$  by their graph invariants concretely.
\end{abstract}
%

%
%
\sftt{Introduction}
\label{section:first}
Super selection sectors play important role
in not only local quantum physics (\cite{DHR,Haag})
but also operator algebras (\cite{Izumi1}).
A sector in the theory of operator algebras is defined
as a unitary equivalence class of a unital C$^{*}$-algebra.
Special sectors are classified by the statistical dimension,
which is an additive, positive integer-valued invariant
of sectors.
we introduced endomorphisms of the Cuntz
algebra $\con$ as follows: Let ${\goth S}_{N,l}$ be the
set of all bijections on the set $\nset{l}$.
For $\sigma\in {\goth S}_{N,l}$,
let $\psi_{\sigma}$ be the endomorphism of $\con$
defined by
%
%
\begin{equation}
\label{eqn:end}
\psi_{\sigma}(s_{i})\equiv 
u_{\sigma}s_{i}\quad(i\edot)
\end{equation}
where
$u_{\sigma}\equiv \sum_{J\in \nset{l}}
s_{\sigma(J)}(s_{J})^{*}$ and
$s_{J}\equiv s_{j_{1}}\cdots s_{j_{l}}$
for $J=(j_{1},\ldots,j_{l})$.
$\psi_{\sigma}$ is called the {\it $l$-th order 
permutative endomorphism} of $\con$ by $\sigma$ ($\S$6-1 \cite{AK2}).
Such endomorphism is concrete and naive but,
it is unknown whether it satisfies ingredients
to define the statistical dimension or not.
For example, we show two endomorphisms $\rho_{1},\rho_{2}$ of $\co{3}$ as follows:
\[
\left\{
\begin{array}{rl}
\rho_{1}(s_{1})\equiv& s_{23,1}
+s_{31,2}+s_{12,3},\\
&\\
\rho_{1}(s_{2})\equiv &s_{32,1}
+s_{13,2}+s_{21,3},\\
&\\
\rho_{1}(s_{3})\equiv &s_{11,1}
+s_{22,2}+s_{33,3},\\
\end{array}
\right.
\quad
\left\{
\begin{array}{rl}
\rho_{2}(s_{1})\equiv& s_{32,1}+s_{11,2}+s_{12,3},\\
&\\
\rho_{2}(s_{2})\equiv &s_{33,1}+s_{22,2}+s_{13,3},\\
&\\
\rho_{2}(s_{3})\equiv &s_{23,1}+s_{21,2}+s_{31,3}\\
\end{array}
\right.
\]
where
$s_{ij,k}\equiv s_{i}s_{j}s_{k}^{*}$ for $i,j,k=1,2,3$.
N.Nakanishi found $\rho_{1}$ by trial and error.
By generalizing $\rho_{1}$, we obtain $\psi_{\sigma}$ in (\ref{eqn:end}).

\noindent
{\bf Question.}
{\it Whether are $\rho_{1}$ and $\rho_{2}$ unitarily equivalent or not?}\\
Branching laws of permutative representations
of $\con$ by permutative endomorphisms 
are computed by directed regular graph (\cite{PE03}).
Since such \brl s are invariant up to unitary equivalence
of endomorphisms, we can classify sectors by their \brl s.
In this article, we introduce invariants of sectors of $\con$
in order to classify sectors more easily than computing \brl s directly.
%
%
\begin{Thm}
\label{Thm:mainzero}
For $\sigma\in {\goth S}_{N,l}$, define a non negative integer-valued 
$N^{l-1}\times N^{l-1}$
matrix $A_{\sigma}=(a_{JK})_{J,K\in\nset{l-1}}$ such that
%
%
\begin{equation}
\label{eqn:maunit}
a_{JK}\equiv \#\{(n,m)\in\nset{2}:\sigma(n,K)=(J,m)\}
\end{equation}
and  non negative real-valued numbers
$c_{1}[\sigma],c_{2}[\sigma],c_{3}[\sigma]$ by
%
%
\begin{equation}
\label{eqn:index}
c_{1}[\sigma]\equiv {\rm Tr}A_{\sigma},\quad
c_{2}[\sigma]\equiv \frac{1}{2}{\rm Tr}(A_{\sigma}^{2}+A_{\sigma}),\quad
c_{3}[\sigma]\equiv \frac{1}{3}{\rm Tr}(A_{\sigma}^{3}+2A_{\sigma}).
\end{equation}
Under these definitions,
if $\psi_{\sigma}\sim\psi_{\eta}$, then
$c_{i}[\sigma]=c_{i}[\eta]$ for $i=1,2,3$.
\end{Thm}

In \cite{PE01},
we classify elements in $\{\psi_{\sigma}:\sigma\in {\goth S}_{2,2}\}$ completely
by computing \brl s for every $\psi_{\sigma}$.
However, it is not good to compute \brl s for every element in 
$\{\psi_{\sigma}:\sigma\in{\goth S}_{3,2}\}$
because $\#{\goth S}_{3,2}=9!\doteq 3.6\times 10^{5}$.
By Theorem \ref{Thm:mainzero}, we obtain a classification of 
this case.
%
%
\begin{Thm}
\label{Thm:maintwo}
For $\sigma\in {\goth S}_{3,2}$ and 
$A_{\sigma}=(a_{ij})_{i,j=1}^{3}$ in Theorem \ref{Thm:mainzero},
define $g_{\sigma}$ by the directed graph with $3$ vertices $v_{1},v_{2},v_{3}$,
$3$ outgoing edges and $3$ incoming edges, 
such that the number of edges from $v_{i}$ to $v_{j}$ 
is $a_{ij}$.
Then the following holds:
\begin{enumerate}
\item
$g_{\sigma}\sim g_{\eta}$ if and only if
$(c_{1}[\sigma],c_{2}[\sigma])=(c_{1}[\eta],c_{2}[\eta])$
where 
$g_{\sigma}\sim g_{\eta}$ means that there is
a permutation matrix $T$ such that
$TA_{\sigma}T^{-1}=A_{\eta}$.
\item
$\#(\{g_{\sigma}:\sigma\in {\goth S}_{3,2}\}/\!\!\!\sim)\,\,=16$.
\item
If $\psi_{\sigma}\sim\psi_{\eta}$, then 
$g_{\sigma}=g_{\eta}$.
\end{enumerate}
\end{Thm}
By using Theorem \ref{Thm:maintwo},
we answer the previous question.
Since $\rho_{1}=\psi_{\sigma}$ and $\rho_{2}=\psi_{\eta}$
for permutations $\sigma,\eta\in {\goth S}_{3,2}$ defined by
\[
\sigma
\left(
\begin{array}{ccc}
11&12&13\\
21&22&23\\
31&32&33\\
\end{array}
\right)
\equiv
\left(
\begin{array}{ccc}
23&31&12\\
32&13&21\\
11&22&33\\
\end{array}
\right),
\quad
\eta
\left(
\begin{array}{ccc}
11&12&13\\
21&22&23\\
31&32&33\\
\end{array}
\right)\equiv
\left(
\begin{array}{ccc}
32&11&12\\
33&22&13\\
23&21&31\\
\end{array}
\right),
\]
graphs $g_{\sigma}$ and $g_{\eta}$ are given as follows:

\noindent
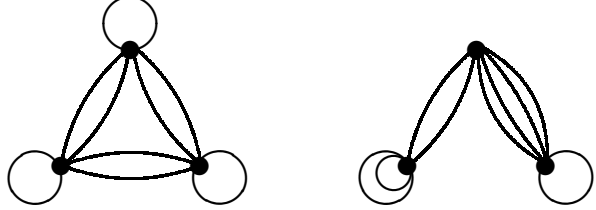
\def\bcircle{\circle*{100}}
\def\numbers{
\put(-400,0){\bcircle}
\put(0,620){\bcircle}
\put(400,0){\bcircle}
}
\def\triplecircle{
\put(0,-30){\circle{250}}
\put(0,-70){\circle{150}}
\put(0,0){\circle{300}}
}
\def\edgethree{
\put(-650,-20){\vector(-2,1){0}}
\put(650,30){\vector(2,-1){0}}
}
\def\edges{
\put(-400,0){\line(1,0){800}}
\put(-400,0){\line(2,3){400}}
\put(400,0){\line(-2,3){400}}
}
\def\sixthloop{
\qbezier(0,0)(400,250)(800,0)
\qbezier(0,0)(400,150)(800,0)
\qbezier(0,0)(400,50)(800,0)
\qbezier(0,0)(400,-50)(800,0)
\qbezier(0,0)(400,-150)(800,0)
\qbezier(0,0)(400,-250)(800,0)
}
\def\fourthloop{
\qbezier(0,0)(400,250)(800,0)
\qbezier(0,0)(400,100)(800,0)
\qbezier(0,0)(400,-100)(800,0)
\qbezier(0,0)(400,-250)(800,0)
}
\def\secondloop{
\qbezier(0,0)(400,250)(800,0)
\qbezier(0,0)(400,-250)(800,0)
}
%
\def\twoone{
\put(-2850,-350){$(3,3,3)$}
\put(-2600,0){\numbers}
\put(-2600,750){\triplecircle}
\put(-3000,120){\triplecircle}
\put(-2200,120){\triplecircle}
}
\def\twotwo{
\put(-2600,0){\numbers}
\put(-2600,750){\triplecircle}
\put(-3000,0){\sixthloop}
\put(-2850,-350){$(3,0,0)$}
}
\def\twothree{
\put(-2600,0){\numbers}
\put(-2600,750){\triplecircle}
\put(-3000,0){\fourthloop}
\put(-3100,0){\circle{250}}
\put(-2100,0){\circle{250}}
\put(-2850,-350){$(3,1,1)$}
}
\def\twofour{
\put(-2600,0){\numbers}
\put(-2600,750){\triplecircle}
\put(-3000,0){\secondloop}
\put(-3100,0){\circle{250}}
\put(-3070,0){\circle{150}}
\put(-2100,0){\circle{250}}
\put(-2130,0){\circle{150}}
\put(-2850,-350){$(3,2,2)$}
}

\def\numbersb{
\put(-400,0){\bcircle}
\put(0,670){\bcircle}
\put(400,0){\bcircle}
}
\def\circlee{\put(0,0){\circle{300}}}
\def\circleeb{
\put(0,0){\circle{300}}
\put(40,150){\vector(1,0){0}}
}
\def\doubleloop{
\qbezier(0,0)(400,150)(800,0)
\qbezier(0,0)(400,-150)(800,0)
}
\def\doubleloopb{
\qbezier(0,0)(400,150)(800,0)
\put(00,0){\line(1,0){800}}
\qbezier(0,0)(400,-150)(800,0)
}
\def\fifthloop{
\qbezier(0,0)(400,250)(800,0)
\qbezier(0,0)(400,150)(800,0)
\put(00,0){\line(1,0){800}}
\qbezier(0,0)(400,-150)(800,0)
\qbezier(0,0)(400,-250)(800,0)
}
\def\doublecircle{
\put(0,-50){\circle{200}}
\put(0,0){\circle{300}}
}
\def\qline{\put(00,0){\line(1,0){800}}}
\def\tripleloop{
\qbezier(0,0)(400,150)(800,0)
\put(00,0){\line(1,0){800}}
\qbezier(0,0)(400,-150)(800,0)
}
\def\edone{
\put(470,0){\vector(1,0){0}}
\put(470,80){\vector(1,0){0}}
\put(470,-80){\vector(1,0){0}}
}
\def\edoneb{
\put(470,80){\vector(1,0){0}}
\put(470,-80){\vector(1,0){0}}
}
\def\edtwo{
\put(350,0){\vector(-1,0){0}}
\put(470,80){\vector(1,0){0}}
\put(470,-80){\vector(1,0){0}}
}
\def\edtwob{
\put(470,80){\vector(1,0){0}}
\put(360,-80){\vector(-1,0){0}}
}
%
%
\def\threefive{
\put(-2600,0){\numbersb}
\put(-3000,0){\doubleloop}
\put(-2260,0){\rotatebox{120}{\doubleloop}}
\put(-2670,690){\rotatebox{-120}{\doubleloop}}
\put(-2600,820){\circlee}
\put(-3150,-70){\rotatebox{120}{\circlee}}
\put(-2220,-70){\rotatebox{-120}{\circlee}}
}
\def\threefiveb{
\put(0,0){\threefive}
\put(-2640,-350){$g_{\sigma}$}
}
\def\threesix{
\put(-2600,0){\numbersb}
\put(-3120,-70){\rotatebox{120}{\doublecircle}}
\put(-2220,-70){\rotatebox{-120}{\circlee}}
\put(-2300,0){\rotatebox{120}{\fourthloop}}
\put(-2670,690){\rotatebox{-120}{\doubleloop}}
\put(-2650,-350){$g_{\eta}$}
}
%
%

\noindent
%
\setlength{\unitlength}{.023013mm}
\begin{picture}(3011,1050)(-3199,-230)
\thicklines
\put(1000,0){\threefiveb}
\put(3000,0){\threesix}
\end{picture}

\noindent
Hence $\rho_{1}\not\sim \rho_{2}$.

In $\S$\ref{section:second},
we show an algorithm to compute \brl s 
of permutative representations by permutative
endomorphisms.
In $\S$\ref{section:third}, we introduce a new index for
more general sectors.
We show formulae among graph invariants and such indices. 
By these formulae, we prove Theorem \ref{Thm:mainzero}.
In $\S$\ref{section:fourth}, we show examples.
We classify $\psi_{\sigma}$ for 
$\sigma\in {\goth S}_{N,l}$ when $(N,l)=(2,2),(2,3),(3,2)$.
We prove Theorem \ref{Thm:maintwo} in $\S$\ref{subsection:fourthtwo}
by illustrating all of 16 graphs.

%
%
\sftt{Automaton computing of \brl s}
\label{section:second}
For $N\geq 2$,
let $\con$ be the {\it Cuntz algebra} (\cite{C}), that is,
the C$^{*}$-algebra which is universally generated by
$s_{1},\ldots,s_{N}$ satisfying $s^{*}_{i} s_j=\delta_{ij}I$ for 
$i,j\edot$ and $s_1 s^{*}_1+\cdots+s_N s^{*}_N=I$.
In this article, any representation and endomorphism
are assumed unital and $*$-preserving.
Two endomorphisms $\rho$ and $\rho^{'}$ are {\it equivalent}
if there is a unitary $u\in\con$
such that $\rho^{'}={\rm Ad}u\circ \rho$.
In this case, we denote $\rho\sim\rho^{'}$.

Let $\nset{*}_{1}\equiv \coprod_{k\geq 1}\nset{k}$,
$\nset{k}\equiv  \{(j_{n})_{n=1}^{k}:j_{n}\edot$,
$n=1,\ldots,k\}$ for $k\geq 1$.
For $J=(j_{l})_{l=1}^{k}\in\nset{*}_{1}$, 
a representation $({\cal H},\pi)$ of $\con$ is 
$P(J)$ if there is a unit cyclic vector $\Omega\in {\cal H}$
such that $\pi(s_{J})\Omega=\Omega$
and $\{\pi(s_{j_{l}}\cdots s_{j_{k}})\Omega\}_{l=1}^{k}$
is an orthogonal family. 
Such representation exists uniquely up to unitary equivalence.
Hence the symbol $P(J)$ makes sense as an equivalence
class of representations.
$P(J)$ is equivalent to a cyclic permutative representation
of $\con$ with a cycle in \cite{BJ,DaPi2}.
When $({\cal H},\pi)$ is $P(J)$, we denote $\pi\circ \psi_{\sigma}$
by $P(J)\circ \psi_{\sigma}$
for $\psi_{\sigma}$ in (\ref{eqn:end}).
In \cite{PE01}, we show that for each $J$,
there are $J_{1},\ldots,J_{m}$,
$1\leq m\leq N^{l-1}$ such that
$P(J)\circ \psi_{\sigma}$ is uniquely decomposed 
into the direct sum of $P(J_{1}),\ldots,P(J_{m})$
up to unitary equivalence:
%
%
\begin{equation}
\label{eqn:branching}
P(J)\circ \psi_{\sigma}\sim
P(J_{1})\oplus \cdots\oplus P(J_{m}).
\end{equation}
Concrete several \brl s by $\psi_{\sigma}$
are given in \cite{PE01}. 

According to \cite{PE03}, 
we show an algorithm to seek $J_{1},\ldots,J_{m}$
in (\ref{eqn:branching})
for $J$ by reducing problem to
a semi-Mealy machine as an input ($=J$)
and outputs ($=J_{1},\ldots,J_{m}$).
A {\it semi-Mealy machine} is a data
$(Q,\Sigma,\Delta,\delta,\lambda)$
which consists of nonempty finite sets
$Q,\Sigma,\Delta$ and two maps
$\delta$ from $Q\times\Sigma^{*}$ to $Q$,
$\lambda$ from $Q\times\Sigma^{*}$ to $\Delta^{*}$
where $\Sigma^{*}$ and $\Delta^{*}$ are free semigroups 
generated by $\Sigma$ and $\Delta$, respectively
and $\delta(q,wa)\equiv \delta(\delta(q,w),a)$
and $\lambda(q,wa)\equiv \lambda(q,w)\lambda(\delta(q,w),a)$
for $q\in Q$, $w\in \Sigma^{*}$ and $a\in \Sigma$.
For symbols $a_{1},\ldots,a_{N}$, $b_{1},\ldots,b_{N}$,
$J=(j_{1},\ldots,j_{k})\in\nset{k}$, $r\geq 1$,
we denote $a_{J}\equiv a_{j_{1}}\cdots a_{j_{k}}$,
$b_{J}\equiv b_{j_{1}}\cdots b_{j_{k}}$
and $a_{J}^{r}\equiv a_{J}\cdots a_{J}$ ($r$-times).
The {\it Mealy diagram} ${\cal D}(M)$ of a semi-Mealy machine
$M=(Q,\Sigma,\Delta,\delta,\lambda)$
is a directed graph with labeled edges which has
a set $Q$ of vertices
and a set $E\equiv \{(q,\delta(q,a),a)\in
Q\times Q\times \Sigma:q\in Q,\, a\in \Sigma\}$
of directed edges with labels.
The meaning of $(q,\delta(q,a),a)$
is an edge from $q$ to $\delta(q,a)$
with a label $a/\lambda(q,a)$ for $a\in\Sigma$:

\noindent
%
%
\def\jei{
  \put(0,0){\oval(600,300)}
  \put(-70,-40){$q$}
}
\def\jeitwo{
 \put(0,0){\oval(600,300)}
 \put(-220,-50){$\delta(q,a)$}
}
\def\figures{
   \put(0,0){\jei}
   \put(1500,0){\jeitwo}
\put(300,0){\vector(1,0){900}}
\put(500,40){$^{a/\lambda(q,a)}$}
}
%
%
\setlength{\unitlength}{.022566mm}
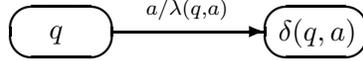
\begin{picture}(1001,280)(199,-100)
\thicklines
\put(2000,0){\figures}
\end{picture}
%

\noindent
A sequence $C=(q_{i_{1}},\ldots,q_{i_{k}})$
in $Q$ is a {\it $k$-cycle} in $M$ by $a_{j_{1}}\cdots a_{j_{k}}\in\Sigma^{*}$ 
if $q_{i_{1}},\ldots,q_{i_{k}}$ satisfy that
$\delta(q_{i_{t}},a_{j_{t}})=q_{i_{t+1}}$ for $t=1,\ldots,k-1$
and $\delta(q_{i_{k}},a_{j_{k}})=q_{i_{1}}$ when $k\geq 2$,
and $\delta(q_{i_{1}},a_{j_{1}})=q_{i_{1}}$ when $k=1$.
We often denote $C$ by $q_{i_{1}}\cdots q_{i_{k}}$ simply.
For a cycle $q_{i_{1}}\cdots q_{i_{k}}$, we {\it do not} assume that
$q_{i_{t}}\ne q_{i_{t^{'}}}$ when $t\ne t^{'}$ in this article.
For $\sigma\in {\goth S}_{N,l}$ with $l\geq 2$ and $J\in\nset{l}$,
we define $\sigma_{1}(J),\ldots,\sigma_{l}(J)\in
\nset{}$ by
$\sigma(J)=(\sigma_{1}(J),\ldots,\sigma_{l}(J))$
and let $\sigma_{n,m}(J)\equiv (\sigma_{n}(J),\ldots,
\sigma_{m}(J))$ for $1\leq n<m\leq l$.
Define $\nset{0}\equiv \{0\}$ for convenience.
%
%
\begin{defi}
\label{defi:mealy}
For $\sigma\in {\goth S}_{N,l}$,
define a data $M_{\sigma}\equiv (Q,\Sigma,\Delta,\delta,\lambda)$
by three sets $Q\equiv\{q_{K}:K\in \nset{l-1}\}$,
$\Sigma\equiv\{a_{j}\}_{j=1}^{N}$,
$\Delta\equiv \{b_{j}\}_{j=1}^{N}$
and two maps  $\delta:Q\times \Sigma^{*}\to Q$, 
$\lambda:Q\times \Sigma^{*}\to \Delta^{*}$,
\[
\delta(q_{K},a_{i})\equiv
\left\{
\begin{array}{ll}
q_{0} \quad\!&(l=1),\\
&\\
q_{(\sigma^{-1})_{2,l}(K,i)}\quad\!&(l\geq 2),\\
\end{array}
\right.\quad \!\!
\lambda(q_{K},a_{i})\equiv
\left\{
\begin{array}{ll}
b_{\sigma^{-1}(i)} \quad \!&(l=1),\\
&\\
b_{(\sigma^{-1})_{1}(K,i)}\quad \!&(l\geq 2)\\
\end{array}
\right.
\]
for $i\edot$ and $K\in\nset{l-1}$.
$M_{\sigma}\equiv (Q,\Sigma,\Delta,\delta,\lambda)$
is called the {\it semi-Mealy machine 
by $\sigma\in {\goth S}_{N,l}$}.
\end{defi}

\noindent
For $\sigma\in {\goth S}_{N,l}$,
the Mealy diagram ${\cal D}(M_{\sigma})$
of $M_{\sigma}$ is a directed regular graph with 
$N^{l-1}$ vertices, $N$ outgoing edges and $N$ incoming edges
when we forget labels of ${\cal D}(M_{\sigma})$.
For a given $J=(j_{i})_{i=1}^{k}\in\nset{k}$,
define $R_{J}\equiv Q_{J}/\!\!\!\sim$
where $Q_{J}\equiv \{q\in Q:\exists n\in{\bf N}\,
s.t.\, \delta(q,x^{n})=q\}$,
$x\equiv a_{j_{1}}\cdots a_{j_{k}}\in \Sigma^{*}$
and $\sim$ is defined as $q\sim q^{'}$
if there is $n\in {\bf N}\cup \{0\}$
such that $\delta(q,x^{n})=q^{'}$.
Then $R_{J}$ is the set of all cycles in $Q$ 
by the input word $x$ and $R_{J}\ne \emptyset$.

%
%
\begin{Thm}
\label{Thm:mainone}
(\cite{PE03})
Let $M_{\sigma}\equiv (Q,\Sigma,\Delta,\delta,\lambda)$
be in Definition \ref{defi:mealy}.
For $J\in\nset{k}$, 
let $p_{1},\ldots,p_{m}\in Q_{J}$ be all representatives of elements in $R_{J}$
and $r_{i}\equiv {\rm min}\{n\in {\bf N}:\delta(p_{i},x^{n})=p_{i}\}$
for $i=1,\ldots,m$.
If $J_{1},\ldots,J_{m}\in\nset{*}_{1}$
are defined by the output word $b_{J_{i}}=\lambda(p_{i},x^{r_{i}})\in \Delta^{*}$
for $i=1,\ldots,m$, then 
$P(J)\circ \psi_{\sigma}\sim P(J_{1})\oplus \cdots\oplus P(J_{m})$.
\end{Thm}

\noindent
In Theorem \ref{Thm:mainone}, if $p_{1}^{'},\ldots,p_{m}^{'}$
are another representatives of $R_{J}$,
then the associated $J_{1}^{'},\ldots,J_{m}^{'}$
satisfy that $P(J_{i}^{'})\sim P(J_{i})$ for each
$i=1,\ldots,m$.

%
%
\sftt{Cyclic index of endomorphism}
\label{section:third}
For $J=(j_{i})_{i=1}^{k}\in \nset{k}$ and $\tau\in {\bf Z}_{k}$,
define $\tau(J)\equiv (j_{\tau(i)})_{i=1}^{k}$.
For $J_{1},J_{2}\in \nset{*}_{1}$,
$J_{1}\sim J_{2}$ if there are $k\geq 1$
and $\tau\in {\bf Z}_{k}$
such that $J_{1},J_{2}\in\nset{k}$
and $\tau(J_{1})=J_{2}$.
$P(J_{1})\sim P(J_{2})$ if and only if $J_{1}\sim J_{2}$.
For $J_{1}=(j_{i})_{i=1}^{k},
J_{2}=(j_{i}^{'})_{i=1}^{k}\in\nset{k}$,
$J_{1}\prec J_{2}$ if $\sum_{l=1}^{k}(j_{l}^{'}-j_{l})N^{k-l}\geq 0$.
$J\in\nset{*}_{1}$ is {\it minimal} if 
$J\prec J^{'}$ for each $J^{'}\in\nset{*}_{1}$
such that $J\sim J^{'}$.
Especially, any element in $\nset{}$
is minimal.

Let ${\rm End}\con$ be the set of all unital $*$-endomorphisms of $\con$.
For $\rho\in {\rm End}\con$, assume that
%
%
\begin{equation}
\label{eqn:assumptions}
\begin{array}{l}
\forall J,J^{'}\in\nset{*}_{1},\,
\exists m(J|\rho|J^{'})\in \{0,1,2,\ldots\}\cup\{\infty\}\,s.t.\\
P(J)\circ \rho=
\bigoplus_{J^{'}\in X_{N,*}}P(J^{'})^{\oplus m(J|\rho|J^{'})}.
\end{array}
\end{equation}
For such $\rho$, define
a non negative integer $c_{k}(\rho)$ by the sum of multiplicities
%
%
\begin{equation}
\label{eqn:cdef}
c_{k}(\rho)\equiv \sum_{J,J^{'}\in X_{N,k}}
m(J|\rho|J^{'})\quad(k\geq 1)
\end{equation}
where $X_{N,k}$ is the set of all minimal elements in $\nset{*}_{1}$.
If $\rho$ satisfies (\ref{eqn:assumptions}),
then $m(J|\rho|J^{'})$'s are uniquely determined
because the l.h.s of the \brl\ in (\ref{eqn:assumptions})
is also a permutative representation.
Hence $c_{k}(\rho)$ is well-defined.
For example,
any permutative endomorphism $\rho$ satisfies (\ref{eqn:assumptions})
and $c_{k}(\rho)<\infty$ for each $k\geq 1$.
We call $c_{k}(\rho)$ by the {\it $k$-th cyclic index} of $\rho$.
Because \brl s are preserved by unitary equivalence, 
if $\rho$ and $\rho^{'}$ satisfy (\ref{eqn:assumptions}) and
$\rho\sim \rho^{'}$, then $c_{k}(\rho)=c_{k}(\rho^{'})$
for each $k\geq 1$.

For $\rho\in {\rm End}\con$,
define $[\rho]\equiv \{\rho^{'}\in {\rm End}\con:
\rho^{'}\sim \rho\}$. $[\rho]$ is the sector of $\con$
with the representative $\rho$.
We see that $c_{k}$ is well-defined on 
${\rm Sect}_{ad}\con\equiv \{[\rho]:\rho\in {\rm End}\con,\,\rho\mbox{ is admissible}\}$,
that is, $c_{k}$ is an invariant of sectors of $\con$.
Define $H_{N}(\con)\equiv \{(\xi_{i})_{i=1}^{N}\in
(\con)^{N}:\xi_{1},\ldots,\xi_{N}\mbox{ satisfy relations }$
of canonical generators of $\con\}$.
For $\rho_{1},\ldots,\rho_{N}\in {\rm End}\con$
and $\xi=(\xi_{i})_{i=1}^{N}\in H_{N}(\con)$,
define
$<\xi|\rho_{1},\ldots,\rho_{N}>\in {\rm End}\con$ by
\[
<\xi|\rho_{1},\ldots,\rho_{N}>\equiv 
{\rm Ad}\xi_{1}\circ \rho_{1}+\cdots 
+{\rm Ad}\xi_{N}\circ \rho_{N}.\]
Then $[<\xi|\rho_{1},\ldots,\rho_{N}>]
= [<\xi^{'}|\rho_{1},\ldots,\rho_{N}>]$
for each $\xi,\xi^{'}\in H_{N}(\con)$.
$[<\xi|\rho_{1},\ldots,\rho_{N}>]
=[<\xi|\rho_{1}^{'},\ldots,\rho_{N}^{'}>]$
if $\rho_{i}\sim \rho^{'}_{i}$ for each $i=1,\ldots,N$.
$[<\xi|\rho_{1},\ldots,\rho_{N}>]
=
[<\xi|\rho_{\sigma(1)},\ldots,\rho_{\sigma(N)}>]$
for each $\sigma\in {\goth S}_{N}$.
In consequence,
the notation
\[[\rho_{1}]+\cdots +[\rho_{N}]\equiv 
[<\xi|\rho_{1},\ldots,\rho_{N}>]\]
is well-defined and 
$[\rho_{1}]+\cdots +[\rho_{N}]=
[\rho_{\sigma(1)}]+\cdots +[\rho_{\sigma(N)}]$
for each $\sigma\in {\goth S}_{N}$.
Further 
$[\rho_{1}]+\cdots +[\rho_{N}]$ satisfies
the associative law (\cite{SE01}).
We simply denote
$[\rho_{1}]+\cdots +[\rho_{N}]$ by
$\rho_{1}+\cdots +\rho_{N}$.
%
%
\begin{prop}
If $\rho_{1},\ldots,\rho_{N}$ satisfy (\ref{eqn:assumptions}),
then $c_{k}(\rho_{1}+\cdots +\rho_{N})=
c_{k}(\rho_{1})+\cdots+c_{k}(\rho_{N})$ for each $k\geq 1$.
\end{prop}
%
%
\pr
Let $\rho\equiv \rho_{1}\oplus\cdots\oplus \rho_{N}$
and $X_{N,*}\equiv \bigcup_{k\geq 1}X_{N,k}$.
Because $P(J)\circ \rho=P(J)\circ \rho_{1}\oplus \cdots\oplus P(J)\circ \rho_{N}$,
$P(J)\circ  \rho=
\bigoplus_{J^{'}\in X_{N,*}}P(J^{'})^{\oplus\sum_{j=1}^{N} m(J|\rho_{j}|J^{'})}$.
From this, $m(J|\rho|J^{'})=\sum_{j=1}^{N} m(J|\rho_{j}|J^{'})$.
This implies the statement.
\qedh

\noindent
In consequence,
$c_{k}$ is an additive invariant of sectors for each $k\geq 1$.

If $\rho$ is the canonical endomorphism of $\con$,
then $c_{k}(\rho)=\#X_{N,k}\cdot N$ for each $k\geq 1$.
If $\rho$ is a permutation of canonical generators of $\con$,
then  $c_{k}(\rho)=\#X_{N,k}$ for each $k\geq 1$.

For a directed finite graph $g$ with the set $V=\{v_{i}\}_{i=1}^{m}$ of vertices,
$A=(a_{ij})$ is the {\it adjacency matrix of $g$} if
$a_{ij}$ is the number of edges from $v_{i}$ to $v_{j}$.
A sequence $(v_{j_{1}},\ldots, v_{j_{k}})$ in $V$ 
is a {\it $k$-cycle} in $g$ if $(v_{j_{1}},v_{j_{2}}),\ldots,(v_{j_{k-1}},v_{j_{k}}),
(v_{j_{k}},v_{j_{1}})$ are directed edges of $g$.
In this article, we do not assume that $v_{j_{i}}\ne v_{j_{i^{'}}}$ when $i\ne i^{'}$.
If $M=(Q,\Sigma,\Delta,\delta,\lambda)$ is a semi-Mealy machine,
then the Mealy diagram ${\cal D}(M)$ gives a directed 
graph by forgetting input/output labels at each edge of ${\cal D}(M)$. 
If $Q=\{q_{j}\}_{j=1}^{m}$, then
the Mealy diagram ${\cal D}(M)$ of $M$ 
gives a graph $g$ with the adjacency matrix
$A=(a_{ij})_{i,j=1}^{m}$ such that
$a_{ij}=\#\{l\in \nset{}:\delta(q_{i},a_{l})=q_{j}\}$
for $i,j=1,\ldots,m$.
For $\sigma\in {\goth S}_{N,l}$, 
the adjacency matrix $A_{\sigma}=(a_{JK})_{J,K\in \nset{l-1}}$ 
of ${\cal D}(M_{\sigma})$ is (\ref{eqn:maunit}) in $\S$\ref{section:first}.
Define $c_{k}(g)$ by the total number of all $k$-cycles in $g$.
For example,
$c_{1}(g)$ is the total number of all $1$-cycles, that is,
the number of edges which starts a vertex $v$ and returns $v$ again. 
$c_{2}(g)$ is the total number of all $2$-cycles, that is,
the number of paths with length $2$ which starts a vertex $v$ and returns $v$ again.
For the adjacency matrix $A$ of a directed graph $g$, the following holds:
%
%
\begin{equation}
\label{eqn:trace}
c_{1}(g)={\rm Tr}A,\quad c_{2}(g)=\frac{1}{2}{\rm Tr}(A+A^{2}),\quad
c_{3}(g)=\frac{1}{3}{\rm Tr}(A^{3}+2A).
\end{equation}
%
%
\begin{prop}
\label{prop:formula}
If $\rho\in E_{N,l}$ and 
$g$ is the Mealy diagram associated with $\rho$,
then $c_{k}(\rho)=c_{k}(g)$ for each $k\geq 1$.
\end{prop}
%
%
\pr
By Theorem \ref{Thm:mainone} and (\ref{eqn:assumptions}),
the statement holds.
\qedh

\noindent
{\it Proof of Theorem \ref{Thm:mainzero}.}
By (\ref{eqn:index}), (\ref{eqn:trace}) 
and Proposition \ref{prop:formula}, 
we see that
$c_{i}[\sigma]=c_{i}(\psi_{\sigma})$ for $i=1,2,3$.
Hence the statement holds.
\qedh

%
%
\sftt{Example}
\label{section:fourth}
For $\psi_{\sigma}$ in (\ref{eqn:end}), define
%
%
\begin{equation}
\label{eqn:edef}
E_{N,l}\equiv \{\psi_{\sigma}\in \endcon:\sigma
\in {\goth S}_{N,l}\}\quad (l\geq 1).
\end{equation}
We show classifications of elements of $E_{2,2}$, $E_{3,2}$ and $E_{2,3}$.
Let ${\goth g}_{N,l}$ be the set of all directed graphs which is isomorphic
to the Mealy diagram associated with $\sigma\in {\goth S}_{N,l}$.
Then ${\goth g}_{N,l}$ coincides with the set of 
all directed regular graphs with $N^{l-1}$ vertices,
$N$ outgoing edges and $N$ incoming edges.
Define ${\goth S}_{N,l}(g)\equiv \{\sigma\in {\goth S}_{N,l}:
{\cal D}(M_{\sigma})\sim g\}$,
$E_{N,l}(g)\equiv \{\rho\in E_{N,l}:
\exists \sigma\in {\goth S}_{N,l}(g)\,s.t.\,\rho=\psi_{\sigma}\}$
and
$SE_{N,l}(g)\equiv E_{N,l}(g)/\!\!\!\sim$ for $g\in {\goth g}_{N,l}$,
In order to classify elements in 
$E_{N,l}$, we list up graphs in ${\goth g}_{N,l}$ 
at each subsection.

For $\sigma\in {\goth S}_{N,2}$,
\[
\psi_{\sigma}(s_{i})=s_{\sigma(i,1)}s_{1}^{*}+\cdots+
s_{\sigma(i,N)}s_{N}^{*}\quad(i\edot).\]
The semi-Mealy machine 
$M_{\sigma}\equiv (Q,\Sigma,\Delta,\delta,\lambda)$
by $\sigma$ is given as
$Q=\{q_{j}\}_{j=1}^{N}$,
$\Sigma=\{a_{j}\}_{j=1}^{N}$, 
$\Delta= \{b_{j}\}_{j=1}^{N}$,
\[
\delta:Q\times \Sigma\to Q;\,
\delta(q_{i},a_{j})=
q_{(\sigma^{-1})_{2}(i,j)},\quad
\lambda:Q\times \Sigma\to \Delta;\,
\lambda(q_{i},a_{j})=
b_{(\sigma^{-1})_{1}(i,j)}
\]
for $i,j\edot$.
%
%
%
\ssft{$E_{2,2}$}
\label{subsection:fourthone}
${\goth g}_{2,2}$ consists of the following $3$ graphs:

\noindent
%
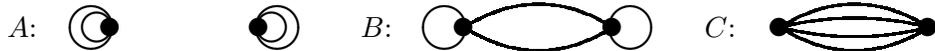
\def\bcircle{\circle*{100}}
\def\numbersw{
\put(-400,0){\bcircle}
\put(400,0){\bcircle}
}
\def\numbers{
\put(-400,0){\bcircle}
\put(400,0){\bcircle}
}
\def\triplecircle{
\put(0,-30){\circle{250}}
\put(0,-70){\circle{150}}
\put(0,0){\circle{300}}
}
\def\sixthloop{
\qbezier(0,0)(400,250)(800,0)
\qbezier(0,0)(400,-250)(800,0)
}
\def\fourloop{
\qbezier(0,0)(400,250)(800,0)
\qbezier(0,0)(400,-250)(800,0)
\qbezier(0,0)(400,100)(800,0)
\qbezier(0,0)(400,-100)(800,0)
}
\def\typeone{
\put(-3550,-50){$B$:}
\put(-2600,0){\numbersw}
\put(-3000,0){\sixthloop}
\put(-3100,0){\circle{250}}
\put(-2100,0){\circle{250}}
}
\def\typetwo{
\put(-3550,-50){$A$:}
\put(-2600,0){\numbersw}
\put(-3100,0){\circle{250}}
\put(-3070,0){\circle{150}}
\put(-2100,0){\circle{250}}
\put(-2130,0){\circle{150}}
}
\def\typethree{
\put(-3400,-50){$C$:}
\put(-2600,0){\numbersw}
\put(-3000,0){\fourloop}
}
%
\setlength{\unitlength}{.024713mm}
\begin{picture}(3011,320)(-3099,-120)
\thicklines
\put(0,0){\typetwo}
\put(1900,0){\typeone}
\put(3600,0){\typethree}
\end{picture}

\noindent
where these directions are unique up to permutation of vertices.
$c_{1}(g)$ for $g=$A, B, C are $4,2,0$, respectively.
The graph of the canonical endomorphisms of $\co{2}$ is A.
The graph of the automorphism of permutation of generators is B.
In this way, these graphs give sketchy classification of 
endomorphisms of $\co{2}$.

In \cite{PE01}, we have more explicit results as follows:
The number of unitary equivalence classes 
of elements in $E_{2,2}$ is $16$.
$G_{2}\equiv {\rm Aut}\co{2}\cap E_{2,2}$ is a subgroup of 
the automorphism group ${\rm Aut}\co{2}$ of $\co{2}$
which is isomorphic to the Klein's four-group.
$G_{2}$ consists of two outer and two inner automorphisms.
$E_{2,2}\setminus G_{2}$ consists of $10$ irreducible 
and $10$ reducible endomorphisms.
The numbers of equivalence classes
are $5$ and $9$, respectively.

Their application for fermion algebra is given by \cite{AK1,AK2}.
These results are given by computing \brl s of every endomorphism concretely. 
Such method is possible because $\#E_{2,2}=24$.
However this is not effective 
when one of $N,k$ is greater than equal $3$ 
because $\#E_{N,k}$ is too large to compute every \brl s. 

%
%
\ssft{$E_{3,2}$}
\label{subsection:fourthtwo}
In order to classify $E_{3,2}$,
we list up Mealy diagrams for every elements in ${\goth S}_{3,2}$.
We see that $M_{\sigma}=(\{q_{1},q_{2},q_{3}\},
\{a_{1},a_{2},a_{3}\}$, $\{b_{1},b_{2},b_{3}\},
\delta,\lambda)$ for each $\sigma\in {\goth S}_{3,2}$.
At the beginning, we show the example $\rho_{1}$ in 
$\S$\ref{section:first}.
By the rule of drawing in $\S$\ref{section:second},
the Mealy diagram is given as follows (\cite{PE01}):

\noindent
%
\def\bcirclex#1{\circle{300}
\put(-50,-20){$q_{#1}$}
}
\def\labee#1#2{
\put(0,0){$a_{#1}/b_{#2}$}
}
\def\numbersx{
\put(-800,0){\bcirclex{2}}
\put(0,1380){\bcirclex{1}}
\put(800,0){\bcirclex{3}}
}
\def\edgethreeq{
\thicklines
\qbezier(-650,-30)(0,-400)(650,-30)
\qbezier(-650,30)(0,400)(650,30)
\put(-650,-20){\vector(-2,1){0}}
\put(650,30){\vector(2,-1){0}}
}
\def\edgesx{
\thicklines
\put(0,0){\edgethreeq}
\put(230,700){\rotatebox{120}{\edgethreeq}}
\put(-560,690){\rotatebox{-120}{\edgethreeq}}
}
\def\semix{
\qbezier(150,0)(0,150)(-150,0)
\qbezier(150,0)(200,-75)(150,-150)
\qbezier(-150,0)(-200,-75)(-150,-150)
\put(-110,-170){\vector(1,-1){0}}}
\def\labelsx{
\put(-180,1750){\labee{1}{3}}
\put(-1340,-300){\labee{2}{3}}
\put(1010,-300){\labee{3}{3}}
%
\put(-1000,750){\labee{3}{1}}
\put(700,750){\labee{2}{1}}
\put(-170,-350){\labee{1}{1}}
\put(-170,30){\labee{1}{2}}
\put(-590,530){\labee{3}{2}}
\put(300,530){\labee{2}{2}}
}
%
\setlength{\unitlength}{.024713mm}
\begin{picture}(3011,2220)(-2099,-330)
\thicklines
\put(0,0){\labelsx}
\put(0,0){\numbersx}
\put(0,0){\edgesx}
%
\put(-20,1620){\semix}
\put(850,-100){\rotatebox{-120}{\semix}}
\put(-1055,-130){\rotatebox{120}{\semix}}
\end{picture}

\noindent
By Theorem \ref{Thm:mainone} and this diagram,
we obtain \brl s by $\rho_{1}$:
\[
\begin{array}{c|c|c|c}
\mbox{input}&\mbox{cycles}&\mbox{outputs}
&\mbox{\brl}\\
\hline
a_{1}&q_{1},q_{2}q_{3}
&b_{3},b_{1}b_{2}&P(1)\circ \rho_{1}=
P(3)\oplus P(12)\\
\hline
a_{1}a_{2}&q_{1}q_{1}q_{3}q_{2}q_{2}q_{3}
&b_{3}b_{1}b_{1}b_{3}b_{2}b_{2}
&P(12)\circ \rho_{1}=
P(113223)\\
\end{array}
\]

For general case, the Mealy diagram is 
a directed regular graphs with $3$ vertices,
$3$ incoming edges and $3$ outgoing edges.
We classify such graphs.
The adjacency matrix $A=(a_{ij})_{i,j=1}^{3}$ is a 
$3\times 3$-matrix with value $0,1,2,3$.
By regularity, $\sum_{i=1}^{3}a_{ij}=\sum_{i=1}^{3}a_{ji}=3$.
Let $a\equiv a_{11},b\equiv a_{22}, c\equiv a_{33}$.
We assume that $3\geq a\geq b\geq c\geq 0$
by permutation of $q_{1},q_{2},q_{3}$.
We see that the possibilities
of $(a,b,c)$ are 
$(3,3,3),(3,2,2),(3,1,1),(3,0,0),(2,2,2)$,
$(2,2,1),(2,1,1)$, $(2,1,0),(2,0,0)$, $(1,1,1)$,
$(1,1,0),(1,0,0),(0,0,0)$.
${\goth g}_{3,2}$ consists of $16$ directed graphs as follows:

\noindent
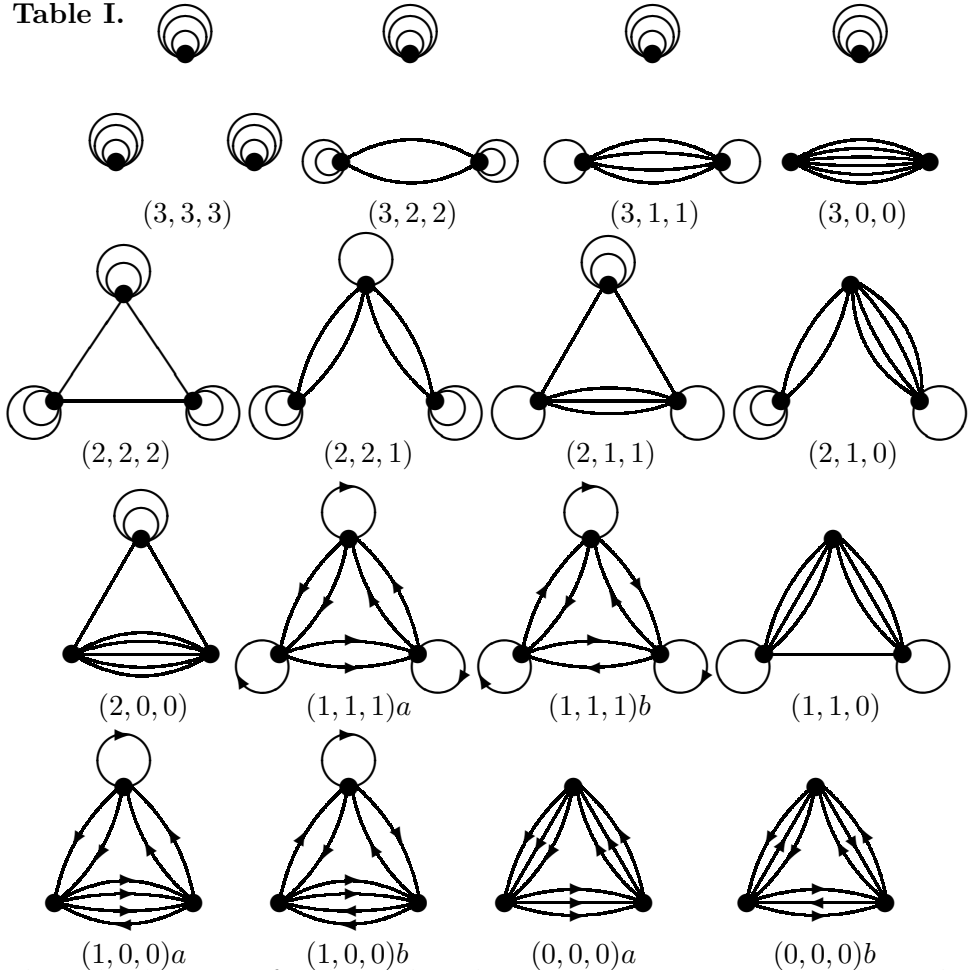
\def\bcircle{\circle*{100}}
\def\numbers{
\put(-400,0){\bcircle}
\put(0,620){\bcircle}
\put(400,0){\bcircle}
}
\def\triplecircle{
\put(0,-30){\circle{250}}
\put(0,-70){\circle{150}}
\put(0,0){\circle{300}}
}
\def\edgethree{
\put(-650,-20){\vector(-2,1){0}}
\put(650,30){\vector(2,-1){0}}
}
\def\edges{
\put(-400,0){\line(1,0){800}}
\put(-400,0){\line(2,3){400}}
\put(400,0){\line(-2,3){400}}
}
\def\sixthloop{
\qbezier(0,0)(400,250)(800,0)
\qbezier(0,0)(400,150)(800,0)
\qbezier(0,0)(400,50)(800,0)
\qbezier(0,0)(400,-50)(800,0)
\qbezier(0,0)(400,-150)(800,0)
\qbezier(0,0)(400,-250)(800,0)
}
\def\fourthloop{
\qbezier(0,0)(400,250)(800,0)
\qbezier(0,0)(400,100)(800,0)
\qbezier(0,0)(400,-100)(800,0)
\qbezier(0,0)(400,-250)(800,0)
}
\def\secondloop{
\qbezier(0,0)(400,250)(800,0)
\qbezier(0,0)(400,-250)(800,0)
}
%
\def\twoone{
\put(-2850,-350){$(3,3,3)$}
\put(-2600,0){\numbers}
\put(-2600,750){\triplecircle}
\put(-3000,120){\triplecircle}
\put(-2200,120){\triplecircle}
}
\def\twotwo{
\put(-2600,0){\numbers}
\put(-2600,750){\triplecircle}
\put(-3000,0){\sixthloop}
\put(-2850,-350){$(3,0,0)$}
}
\def\twothree{
\put(-2600,0){\numbers}
\put(-2600,750){\triplecircle}
\put(-3000,0){\fourthloop}
\put(-3100,0){\circle{250}}
\put(-2100,0){\circle{250}}
\put(-2850,-350){$(3,1,1)$}
}
\def\twofour{
\put(-2600,0){\numbers}
\put(-2600,750){\triplecircle}
\put(-3000,0){\secondloop}
\put(-3100,0){\circle{250}}
\put(-3070,0){\circle{150}}
\put(-2100,0){\circle{250}}
\put(-2130,0){\circle{150}}
\put(-2850,-350){$(3,2,2)$}
}

\def\numbersb{
\put(-400,0){\bcircle}
\put(0,670){\bcircle}
\put(400,0){\bcircle}
}
\def\circlee{\put(0,0){\circle{300}}}
\def\circleeb{
\put(0,0){\circle{300}}
\put(40,150){\vector(1,0){0}}
}
\def\doubleloop{
\qbezier(0,0)(400,150)(800,0)
\qbezier(0,0)(400,-150)(800,0)
}
\def\doubleloopb{
\qbezier(0,0)(400,150)(800,0)
\put(00,0){\line(1,0){800}}
\qbezier(0,0)(400,-150)(800,0)
}
\def\fifthloop{
\qbezier(0,0)(400,250)(800,0)
\qbezier(0,0)(400,150)(800,0)
\put(00,0){\line(1,0){800}}
\qbezier(0,0)(400,-150)(800,0)
\qbezier(0,0)(400,-250)(800,0)
}
\def\doublecircle{
\put(0,-50){\circle{200}}
\put(0,0){\circle{300}}
}
\def\qline{\put(00,0){\line(1,0){800}}}
\def\tripleloop{
\qbezier(0,0)(400,150)(800,0)
\put(00,0){\line(1,0){800}}
\qbezier(0,0)(400,-150)(800,0)
}
\def\edone{
\put(470,0){\vector(1,0){0}}
\put(470,80){\vector(1,0){0}}
\put(470,-80){\vector(1,0){0}}
}
\def\edoneb{
\put(470,80){\vector(1,0){0}}
\put(470,-80){\vector(1,0){0}}
}
\def\edtwo{
\put(350,0){\vector(-1,0){0}}
\put(470,80){\vector(1,0){0}}
\put(470,-80){\vector(1,0){0}}
}
\def\edtwob{
\put(470,80){\vector(1,0){0}}
\put(360,-80){\vector(-1,0){0}}
}
\def\edthree{
\put(0,0){\fourthloop}
\put(470,130){\vector(1,0){0}}
\put(470,50){\vector(1,0){0}}
\put(470,-50){\vector(1,0){0}}
\put(350,-130){\vector(-1,0){0}}
}
\def\edfour{
\put(0,0){\fourthloop}
\put(470,130){\vector(1,0){0}}
\put(470,50){\vector(1,0){0}}
\put(350,-50){\vector(-1,0){0}}
\put(350,-130){\vector(-1,0){0}}
}
%
\def\threeone{
\put(-2600,0){\numbersb}
\put(-3000,0){\doubleloopb}
\put(-2260,0){\rotatebox{120}{\doubleloopb}}
\put(-2670,690){\rotatebox{-120}{\doubleloopb}}
}
\def\threeonea{
\put(0,0){\threeone}
\put(-3000,0){\edone}
\put(-2380,-10){\rotatebox{120}{\edone}}
\put(-2670,690){\rotatebox{-120}{\edone}}
\put(-2850,-350){$(0,0,0)a$}
}
\def\threeoneb{
\put(0,0){\threeone}
\put(-3000,0){\edtwo}
\put(-2380,-10){\rotatebox{120}{\edtwo}}
\put(-2670,690){\rotatebox{-120}{\edtwo}}
\put(-2850,-350){$(0,0,0)b$}
}
\def\threetwo{
\put(-2600,0){\numbersb}
\put(-3000,0){\fourthloop}
\put(-2260,0){\rotatebox{120}{\doubleloop}}
\put(-2670,690){\rotatebox{-120}{\doubleloop}}
\put(-2600,820){\circleeb}
}
\def\threetwoa{
\put(0,0){\threetwo}
\put(-3000,0){\edthree}
\put(-2380,-10){\rotatebox{120}{\edoneb}}
\put(-2670,690){\rotatebox{-120}{\edoneb}}
\put(-2850,-350){$(1,0,0)a$}
}
\def\threetwob{
\put(0,0){\threetwo}
\put(-3000,0){\edfour}
\put(-2380,-10){\rotatebox{120}{\edtwob}}
\put(-2670,690){\rotatebox{-120}{\edtwob}}
\put(-2850,-350){$(1,0,0)b$}
}
\def\threethree{
\put(-2600,0){\numbersb}
\put(-3000,0){\fifthloop}
\put(-2200,0){\rotatebox{120}{\qline}}
\put(-2620,690){\rotatebox{-120}{\qline}}
\put(-2600,800){\doublecircle}
\put(-2850,-350){$(2,0,0)$}
}
\def\threefour{
\put(-2600,0){\numbersb}
\put(-3000,0){\line(1,0){800}}
\put(-2260,0){\rotatebox{120}{\doubleloopb}}
\put(-2670,690){\rotatebox{-120}{\doubleloopb}}
\put(-3120,-70){\rotatebox{120}{\circlee}}
\put(-2220,-70){\rotatebox{-120}{\circlee}}
\put(-2850,-350){$(1,1,0)$}
}
\def\threefive{
\put(-2600,0){\numbersb}
\put(-3000,0){\doubleloop}
\put(-2260,0){\rotatebox{120}{\doubleloop}}
\put(-2670,690){\rotatebox{-120}{\doubleloop}}
\put(-2600,820){\circleeb}
\put(-3360,-70){\rotatebox{120}{\circleeb}}
\put(-2220,-70){\rotatebox{-120}{\circleeb}}
}
\def\threefivea{
\put(0,0){\threefive}
\put(-3000,0){\edoneb}
\put(-2390,0){\rotatebox{120}{\edoneb}}
\put(-2670,690){\rotatebox{-120}{\edoneb}}
\put(-2850,-350){$(1,1,1)a$}
}
\def\threefiveb{
\put(0,0){\threefive}
\put(-3000,0){\edtwob}
\put(-2390,0){\rotatebox{120}{\edtwob}}
\put(-2670,690){\rotatebox{-120}{\edtwob}}
\put(-2850,-350){$(1,1,1)b$}
}
\def\threesix{
\put(-2600,0){\numbersb}
\put(-3120,-70){\rotatebox{120}{\doublecircle}}
\put(-2220,-70){\rotatebox{-120}{\circlee}}
\put(-2300,0){\rotatebox{120}{\fourthloop}}
\put(-2670,690){\rotatebox{-120}{\doubleloop}}
\put(-2850,-350){$(2,1,0)$}
}
\def\threeseven{
\put(-2600,0){\numbersb}
\put(-3000,0){\tripleloop}
\put(-2200,0){\rotatebox{120}{\qline}}
\put(-2620,690){\rotatebox{-120}{\qline}}
\put(-2600,800){\doublecircle}
\put(-3120,-70){\rotatebox{120}{\circlee}}
\put(-2220,-70){\rotatebox{-120}{\circlee}}
\put(-2850,-350){$(2,1,1)$}
}
\def\threeeight{
\put(-2600,0){\numbersb}
\put(-3120,-70){\rotatebox{120}{\doublecircle}}
\put(-2220,-70){\rotatebox{-120}{\doublecircle}}
\put(-2280,0){\rotatebox{120}{\doubleloop}}
\put(-2670,690){\rotatebox{-120}{\doubleloop}}
\put(-2600,820){\circlee}
\put(-2850,-350){$(2,2,1)$}
}
\def\threenine{
\put(-2600,0){\numbers}
\put(-2600,0){\edges}
\put(-2600,750){\doublecircle}
\put(-3120,-70){\rotatebox{120}{\doublecircle}}
\put(-2220,-70){\rotatebox{-120}{\doublecircle}}
\put(-2850,-350){$(2,2,2)$}
}
%
%

\noindent
%
\setlength{\unitlength}{.023013mm}
\begin{picture}(3011,1220)(-3599,-400)
\thicklines
\put(-3600,800){{\bf Table I.}}
\put(0,0){\twoone}
\put(1300,0){\twofour}
\put(2700,0){\twothree}
\put(3900,0){\twotwo}
\end{picture}

\noindent
\begin{picture}(3011,1320)(-3299,-350)
\thicklines
\put(0,0){\threenine}
\put(1400,0){\threeeight}
\put(2800,0){\threeseven}
\put(4200,0){\threesix}
\end{picture}

\noindent
\begin{picture}(3011,1400)(-3399,-300)
\thicklines
\put(0,0){\threethree}
\put(1200,0){\threefivea}
\put(2600,0){\threefiveb}
\put(4000,0){\threefour}
\end{picture}

\noindent
\begin{picture}(3011,1420)(-3299,-300)
\thicklines
\put(0,0){\threetwoa}
\put(1300,0){\threetwob}
\put(2600,0){\threeonea}
\put(4000,0){\threeoneb}
\end{picture}

\noindent
where the direction of each graph without arrows is unique. 
For each graph in Table I,
there are $216$ kinds of configuration of output labels on edges.
For the Mealy diagrams associated with 
$E_{3,2}$, we have the following:\\
\[
\begin{array}{lll}
\begin{array}{l}
\mbox{{\bf Table II.}}\hfill \quad\\
\\
\\
\\
\\
\\
\\
\\
\\
\\
\\
\\
\\
\\
\\
\\
\\
\end{array}
&
\begin{array}{c|c|c|c}
g & c_{1}(g)&c_{2}(g)& \#E_{3,2}(g)/216\\
\hline
(3,3,3) &9 &18 &1 \\
(3,2,2) &7 &13 &27 \\
(3,1,1) &5 &12 &27\\
(3,0,0) &3 &15 &3\\
\hline
(2,2,2) & 6&9& 54\\
(2,2,1) & 5&9 &162 \\
(2,1,1) & 4&7& 324\\
(2,1,0) & 3&9 & 162\\
(2,0,0) & 2&9&54\\
\hline
(1,1,1)a & 3&3&54\\
(1,1,1)b & 3&6&216\\
(1,1,0)  & 2&6&324\\
(1,0,0)a & 1&4&54\\
(1,0,0)b & 1&7&162\\
(0,0,0)a & 0&0&2\\
(0,0,0)b & 0&6&54\\
\end{array}
&\quad\hfill\quad
\end{array}
\]

\noindent
where the number of diagrams is computed as the case may be.

\noindent
{\it Proof of Theorem \ref{Thm:maintwo}.}
By Table II, any element in $E_{3,2}$ belongs to
$E_{3,2}(g)$ for some $g$ appearing in Table I because
we see that $\sum_{g}\#E_{3,2}(g)=9!=\#E_{3,2}$.
By this, Table II and Theorem \ref{Thm:mainzero},
the statement holds. \qedh

\noindent
We show that how the number of  $E_{3,2}(2,1,1)\div 216$ is computed.
The permutation of vertices gives $6$ different diagrams.
Assume that $q_{1}$ has double $1$-cycles and 
there is a directed edge from $q_{1}$ to $q_{2}$.
The choice of the (input) label of the vertex from $q_{1}$ to $q_{2}$ is $3$. 
The choice of label of $1$-cycle at $q_{2}$  is $3$.
The choice of labels of edges from $q_{3}$ is $6$.
Hence the total number of possibilities is given as
$6\times 3\times 3\times  6=324$.

By Table II, we see that $\#E_{3,2}(3,3,3)=216$.
Any element in $E_{3,2}(3,3,3)$ is given as 
$\rho(x)\equiv s_{1}\alpha_{1}(x)s_{1}^{*}+
s_{2}\alpha_{2}(x)s_{2}^{*}+s_{3}\alpha_{3}(x)s_{3}^{*}$
for $x\in \co{3}$
where $\alpha_{1},\alpha_{2},\alpha_{3}\in {\rm Aut}\co{3}$
are permutations of canonical generators.
From this, we can verify that
$\#E_{3,2}(3,3,3)=\{\#{\goth S}_{3}\}^{3}$ again.
The number of equivalence classes in $E_{3,2}(3,3,3)$ is $56$.

%
In consequence, the decomposition 
$SE_{3,2}=\coprod_{g\in {\goth g}_{3,2}}SE_{3,2}(g)$
holds as the disjoint union. Sectors are classified by graphs.
%

%
%
\ssft{$E_{2,3}$}
\label{subsection:fourththree}
For $\sigma\in {\goth S}_{2,3}$,
$M_{\sigma}=(\{q_{11},q_{12},q_{21},q_{22}\},
\{a_{1},a_{2}\},\{b_{1},b_{2}\},\delta,\lambda)$,
$\delta(q_{J},a_{i})=q_{(\sigma^{-1})_{2,3}(J,i)}$,
$\lambda(q_{J},a_{i})=b_{(\sigma^{-1})_{1}(J,i)}$
for $J\in\{1,2\}^{2},\,i=1,2$.
The graph is a regular directed graph with $4$ vertices,
$2$ outgoing edges and $2$ incoming edges.
${\goth g}_{2,3}$ consists of $25$ directed graphs which are classified by
the diagonal part of their adjacency matrices as follows:

\noindent
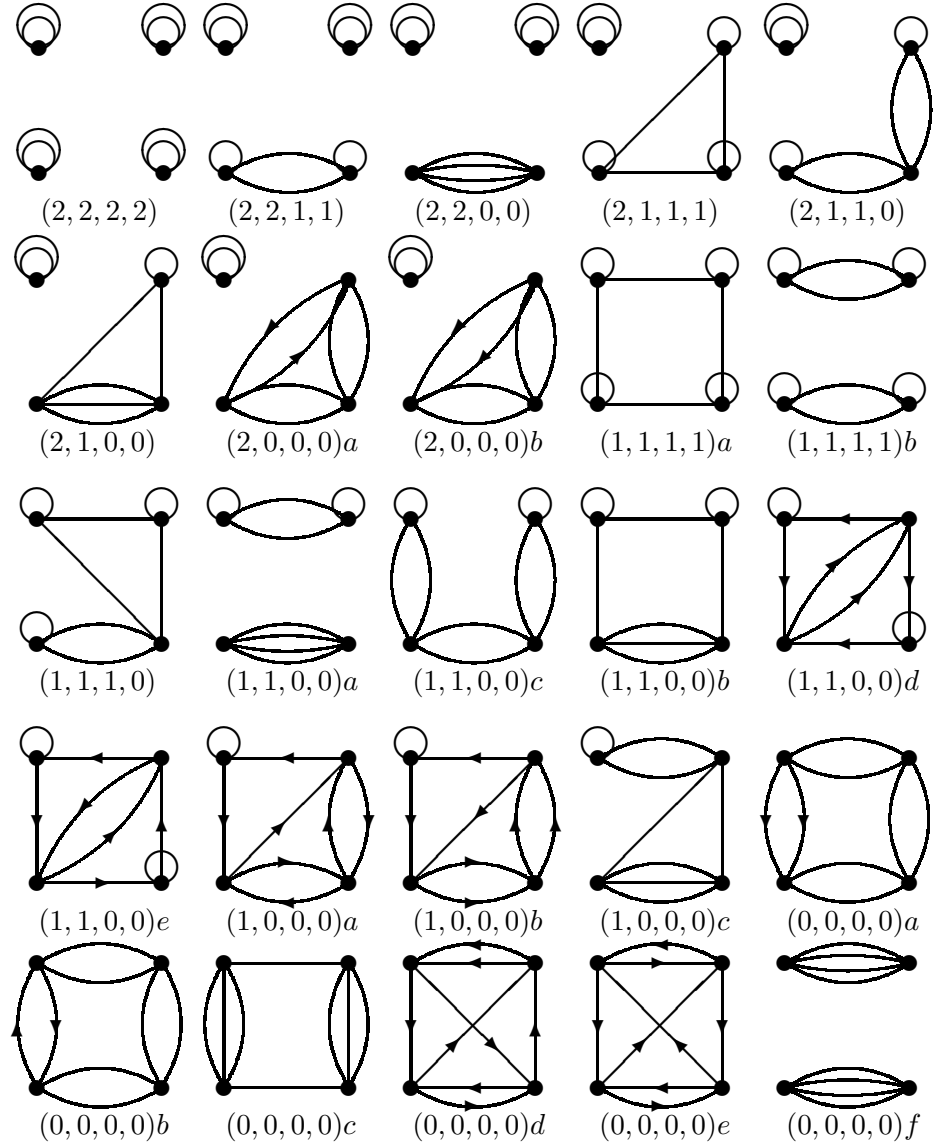
\def\ocircle{\circle{200}}
\def\dcircle{
\put(0,0){\circle{200}}
\put(0,50){\circle{300}}
}
\def\bcircle{\circle*{100}}
\def\numbers{
\put(-400,0){\bcircle}
\put(0,620){\bcircle}
\put(400,0){\bcircle}
}
\def\fourp{
\put(-400,0){\bcircle}
\put(400,0){\bcircle}
\put(-400,800){\bcircle}
\put(400,800){\bcircle}
}
\def\numbersb{
\put(-400,0){\bcircle}
\put(0,670){\bcircle}
\put(400,0){\bcircle}
}
\def\triplecircle{
\put(0,-30){\circle{250}}
\put(0,-70){\circle{150}}
\put(0,0){\circle{300}}
}
\def\edgethree{
\put(-650,-20){\vector(-2,1){0}}
\put(650,30){\vector(2,-1){0}}
}
\def\edges{
\put(-400,0){\line(1,0){800}}
\put(-400,0){\line(2,3){400}}
\put(400,0){\line(-2,3){400}}
}
\def\fouredges{
\put(-400,0){\line(1,0){800}}
\put(-400,0){\line(0,1){800}}
\put(-400,800){\line(1,0){800}}
\put(400,0){\line(0,1){800}}
}
\def\sixthloop{
\qbezier(0,0)(400,250)(800,0)
\qbezier(0,0)(400,150)(800,0)
\qbezier(0,0)(400,50)(800,0)
\qbezier(0,0)(400,-50)(800,0)
\qbezier(0,0)(400,-150)(800,0)
\qbezier(0,0)(400,-250)(800,0)
}
\def\fourthloop{
\qbezier(0,0)(400,250)(800,0)
\qbezier(0,0)(400,100)(800,0)
\qbezier(0,0)(400,-100)(800,0)
\qbezier(0,0)(400,-250)(800,0)
}
\def\secondloop{
\qbezier(0,0)(400,250)(800,0)
\qbezier(0,0)(400,-250)(800,0)
}
\def\secondloopb{
\qbezier(0,0)(580,250)(1160,0)
\qbezier(0,0)(580,-250)(1160,0)
}
\def\curveu{
\qbezier(0,0)(400,250)(800,0)
}
\def\curved{
\qbezier(0,0)(400,-250)(800,0)
}
\def\tn{390}
\def\sgtwotwotwo{
\put(0,0){\fourp}
\put(-\tn,-300){$(2,2,2,2)$}
\put(-400,100){\dcircle}
\put(400,100){\dcircle}
\put(-400,900){\dcircle}
\put(400,900){\dcircle}
}
%
\def\sgtwotwoone{
\put(0,0){\fourp}
\put(-\tn,-300){$(2,2,1,1)$}
\put(-400,900){\dcircle}
\put(400,900){\dcircle}
\put(-400,0){\secondloop}
\put(-400,100){\ocircle}
\put(400,100){\ocircle}
}
%
\def\sgtwotwozero{
\put(0,0){\fourp}
\put(-\tn,-300){$(2,2,0,0)$}
\put(-400,900){\dcircle}
\put(400,900){\dcircle}
\put(-400,0){\fourthloop}
}
%
\def\sgtwooneone{
\put(0,0){\fourp}
\put(-\tn,-300){$(2,1,1,1)$}
\put(-400,900){\dcircle}
\put(400,900){\ocircle}
\put(-400,100){\ocircle}
\put(400,100){\ocircle}
\put(-400,0){\line(1,1){800}}
\put(-400,0){\line(1,0){800}}
\put(400,0){\line(0,1){800}}
}
%
\def\sgtwoonezeroa{
\put(0,0){\fourp}
\put(-\tn,-300){$(2,1,1,0)$}
\put(-400,900){\dcircle}
\put(400,900){\ocircle}
\put(-400,100){\ocircle}
\put(-400,0){\secondloop}
\put(270,0){\rotatebox{90}{\secondloop}}
}
%
\def\sgtwoonezerob{
\put(0,0){\fourp}
\put(-\tn,-300){$(2,1,0,0)$}
\put(-400,900){\dcircle}
\put(400,900){\ocircle}
\put(-400,0){\line(1,1){800}}
\put(-400,0){\line(1,0){800}}
\put(400,0){\line(0,1){800}}
\put(-400,0){\secondloop}
}
%
\def\sgtwozerozeroa{
\put(0,0){\fourp}
\put(-\tn,-300){$(2,0,0,0)a$}
\put(-400,900){\dcircle}
\put(-400,0){\secondloop}
\put(-140,450){\vector(-1,-1){0}}
\put(120,350){\vector(1,1){0}}
\put(270,0){\rotatebox{90}{\secondloop}}
\put(-500,0){\rotatebox{45}{\secondloopb}}
}
\def\sgtwozerozerob{
\put(0,0){\fourp}
\put(-140,450){\vector(-1,-1){0}}
\put(20,250){\vector(-1,-1){0}}
\put(-\tn,-300){$(2,0,0,0)b$}
\put(-400,900){\dcircle}
\put(-400,0){\secondloop}
\put(270,0){\rotatebox{90}{\secondloop}}
\put(-500,0){\rotatebox{45}{\secondloopb}}
}
%
\def\sgoneoneonea{
\put(0,0){\fourp}
\put(-\tn,-300){$(1,1,1,1)a$}
\put(-400,900){\ocircle}
\put(400,900){\ocircle}
\put(-400,100){\ocircle}
\put(400,100){\ocircle}
\put(-400,0){\line(0,1){800}}
\put(-400,0){\line(1,0){800}}
\put(400,0){\line(0,1){800}}
\put(-400,800){\line(1,0){800}}
}
%
\def\sgoneoneoneb{
\put(0,0){\fourp}
\put(-\tn,-300){$(1,1,1,1)b$}
\put(-400,900){\ocircle}
\put(400,900){\ocircle}
\put(-400,100){\ocircle}
\put(400,100){\ocircle}
\put(-400,0){\secondloop}
\put(-400,800){\secondloop}
}
%
\def\sgoneonezero{
\put(0,0){\fourp}
\put(-\tn,-300){$(1,1,1,0)$}
\put(-400,900){\ocircle}
\put(400,900){\ocircle}
\put(-400,100){\ocircle}
\put(-400,0){\secondloop}
\put(-400,800){\line(1,0){800}}
\put(400,0){\line(-1,1){800}}
\put(400,0){\line(0,1){800}}
}
%
\def\sgonezerozeroa{
\put(0,0){\fourp}
\put(-\tn,-300){$(1,1,0,0)a$}
\put(-400,900){\ocircle}
\put(400,900){\ocircle}
\put(-400,0){\fourthloop}
\put(-400,800){\secondloop}
}
%
\def\sgonezerozerob{
\put(0,0){\fourp}
\put(-\tn,-300){$(1,1,0,0)b$}
\put(-400,900){\ocircle}
\put(400,900){\ocircle}
\put(-400,0){\secondloop}
\put(-400,800){\line(1,0){800}}
\put(-400,0){\line(1,0){800}}
\put(-400,0){\line(0,1){800}}
\put(400,0){\line(0,1){800}}
}
\def\sgonezerozeroc{
\put(0,0){\fourp}
\put(-\tn,-300){$(1,1,0,0)c$}
\put(-400,900){\ocircle}
\put(400,900){\ocircle}
\put(-400,0){\secondloop}
\put(270,0){\rotatebox{90}{\secondloop}}
\put(-530,0){\rotatebox{90}{\secondloop}}
}
\def\sgonezerozerod{
\put(-40,550){\vector(1,1){0}}
\put(120,350){\vector(1,1){0}}
\put(-80,800){\vector(-1,0){0}}
\put(-400,330){\vector(0,-1){0}}
\put(-80,00){\vector(-1,0){0}}
\put(400,330){\vector(0,-1){0}}
\put(0,0){\fourp}
\put(-\tn,-300){$(1,1,0,0)d$}
\put(-400,900){\ocircle}
\put(400,100){\ocircle}
\put(-400,0){\line(0,1){800}}
\put(-400,0){\line(1,0){800}}
\put(-400,800){\line(1,0){800}}
\put(400,0){\line(0,1){800}}
\put(-500,0){\rotatebox{45}{\secondloopb}}
}
\def\sgonezerozeroe{
\put(-140,450){\vector(-1,-1){0}}
\put(120,350){\vector(1,1){0}}
\put(-80,800){\vector(-1,0){0}}
\put(-400,330){\vector(0,-1){0}}
\put(80,00){\vector(1,0){0}}
\put(400,430){\vector(0,1){0}}
\put(0,0){\fourp}
\put(-\tn,-300){$(1,1,0,0)e$}
\put(-400,900){\ocircle}
\put(400,100){\ocircle}
\put(-400,0){\line(0,1){800}}
\put(-400,0){\line(1,0){800}}
\put(-400,800){\line(1,0){800}}
\put(400,0){\line(0,1){800}}
\put(-500,0){\rotatebox{45}{\secondloopb}}
}
\def\sgonezerozerota{
\put(-400,330){\vector(0,-1){0}}
\put(-50,800){\vector(-1,0){0}}
\put(0,400){\vector(1,1){0}}
\put(50,130){\vector(1,0){0}}
\put(-50,-130){\vector(-1,0){0}}
\put(270,430){\vector(0,1){0}}
\put(530,330){\vector(0,-1){0}}
\put(0,0){\fourp}
\put(-\tn,-300){$(1,0,0,0)a$}
\put(-400,900){\ocircle}
\put(-400,0){\secondloop}
\put(-400,0){\line(1,1){800}}
\put(-400,0){\line(0,1){800}}
\put(-400,800){\line(1,0){800}}
\put(270,0){\rotatebox{90}{\secondloop}}
}
\def\sgonezerozerotb{
\put(-400,330){\vector(0,-1){0}}
\put(-50,800){\vector(-1,0){0}}
\put(0,400){\vector(-1,-1){0}}
\put(50,130){\vector(1,0){0}}
\put(50,-130){\vector(1,0){0}}
\put(270,430){\vector(0,1){0}}
\put(530,430){\vector(0,1){0}}
\put(0,0){\fourp}
\put(-\tn,-300){$(1,0,0,0)b$}
\put(-400,900){\ocircle}
\put(-400,0){\secondloop}
\put(-400,0){\line(1,1){800}}
\put(-400,0){\line(0,1){800}}
\put(-400,800){\line(1,0){800}}
\put(270,0){\rotatebox{90}{\secondloop}}
}
\def\sgonezerozerotc{
\put(0,0){\fourp}
\put(-\tn,-300){$(1,0,0,0)c$}
\put(-400,900){\ocircle}
\put(-400,800){\secondloop}
\put(-400,0){\line(1,1){800}}
\put(400,0){\line(0,1){800}}
\put(-400,0){\line(1,0){800}}
\put(-400,0){\secondloop}
}
%
\def\sgonea{
\put(-280,330){\vector(0,-1){0}}
\put(-520,330){\vector(0,-1){0}}
\put(0,0){\fourp}
\put(-\tn,-300){$(0,0,0,0)a$}
\put(-400,0){\secondloop}
\put(-400,800){\secondloop}
\put(-530,0){\rotatebox{90}{\secondloop}}
\put(270,0){\rotatebox{90}{\secondloop}}
}
%
\def\sgoneb{
\put(-280,330){\vector(0,-1){0}}
\put(-530,430){\vector(0,1){0}}
\put(0,0){\fourp}
\put(-\tn,-300){$(0,0,0,0)b$}
\put(-400,0){\secondloop}
\put(-400,800){\secondloop}
\put(-530,0){\rotatebox{90}{\secondloop}}
\put(270,0){\rotatebox{90}{\secondloop}}
}
%
\def\sgtwo{
\put(0,0){\fourp}
\put(0,0){\fouredges}
\put(-\tn,-300){$(0,0,0,0)c$}
\put(-530,0){\rotatebox{90}{\secondloop}}
\put(270,0){\rotatebox{90}{\secondloop}}
}
%
\def\sgzerozerofive{
\put(170,230){\vector(1,-1){0}}
\put(-400,330){\vector(0,-1){0}}
\put(-50,800){\vector(-1,0){0}}
\put(-50,920){\vector(-1,0){0}}
\put(-80,320){\vector(1,1){0}}
\put(400,430){\vector(0,1){0}}
\put(-50,0){\vector(-1,0){0}}
\put(50,-120){\vector(1,0){0}}
\put(0,0){\fourp}
\put(0,0){\fouredges}
\put(-400,0){\line(1,1){800}}
\put(400,0){\line(-1,1){800}}
\put(-400,800){\curveu}
\put(-400,0){\curved}
\put(-\tn,-300){$(0,0,0,0)d$}
}
%
\def\sgzerozerofiveb{
\put(80,320){\vector(-1,1){0}}
\put(-400,330){\vector(0,-1){0}}
\put(50,800){\vector(1,0){0}}
\put(-50,920){\vector(-1,0){0}}
\put(-80,320){\vector(1,1){0}}
\put(400,330){\vector(0,-1){0}}
\put(-50,0){\vector(-1,0){0}}
\put(50,-120){\vector(1,0){0}}
\put(0,0){\fourp}
\put(0,0){\fouredges}
\put(-400,0){\line(1,1){800}}
\put(400,0){\line(-1,1){800}}
\put(-400,800){\curveu}
\put(-400,0){\curved}
\put(-\tn,-300){$(0,0,0,0)e$}
}
\def\sgthree{
\put(0,0){\fourp}
\put(-\tn,-300){$(0,0,0,0)f$}
\put(-400,0){\fourthloop}
\put(-400,800){\fourthloop}
}
%
\def\move{600}

\noindent
\setlength{\unitlength}{.020713mm}
\begin{picture}(3011,1720)(-550,-500)
\thicklines
\put(-600,1200){{\bf Table III.}}
\put(0,0){\sgtwotwotwo}
\put(1200,0){\sgtwotwoone}
\put(2400,0){\sgtwotwozero}
\put(3600,0){\sgtwooneone}
\put(4800,0){\sgtwoonezeroa}
\end{picture}

\noindent
\begin{picture}(3011,1520)(-\move,-550)
\thicklines
\put(0,0){\sgtwoonezerob}
\put(1200,0){\sgtwozerozeroa}
\put(2400,0){\sgtwozerozerob}
\put(3600,0){\sgoneoneonea}
\put(4800,0){\sgoneoneoneb}
\end{picture}

\noindent
\begin{picture}(3011,1520)(-\move,-550)
\thicklines
\put(0,0){\sgoneonezero}
\put(1200,0){\sgonezerozeroa}
\put(2400,0){\sgonezerozeroc}
\put(3600,0){\sgonezerozerob}
\put(4800,0){\sgonezerozerod}
\end{picture}

\noindent
\begin{picture}(3011,1520)(-\move,-550)
\thicklines
\put(0,0){\sgonezerozeroe}
\put(1200,0){\sgonezerozerota}
\put(2400,0){\sgonezerozerotb}
\put(3600,0){\sgonezerozerotc}
\put(4800,0){\sgonea}
\end{picture}

\noindent
\begin{picture}(3011,1000)(-\move,-250)
\thicklines
\put(0,0){\sgoneb}
\put(1200,0){\sgtwo}
\put(2400,0){\sgzerozerofive}
\put(3600,0){\sgzerozerofiveb}
\put(4800,0){\sgthree}
\end{picture}

\noindent 
where the direction of each graph without arrows is unique. 
For each graph in Table III,
there are $16$ kinds of configuration of output labels on edges.
Along with Table II, we have the following: 
\vspace{-3mm}
\[\begin{array}{lcl}
\begin{array}{l}
\mbox{{\bf Table IV.}\quad\quad}\\
\\
\\
\\
\\
\\
\\
\\
\\
\\
\\
\\
\\
\\
\\
\\
\\
\\
\\
\\
\\
\\
\\
\\
\\
\\
\end{array}
&
\begin{array}{c|c|c|c}
g & c_{1}(g)&c_{2}(g)& \#E_{2,3}(g)/16\\
\hline
(2,2,2,2)    &8 &12 &1 \\
(2,2,1,1)    &6 &9  &24 \\
(2,2,0,0)    &4 &10 &6 \\
(2,1,1,1)    &5 &6  &64 \\
(2,1,1,0)    &4 &7  &96 \\
(2,1,0,0)    &3 &6 &96 \\
\hline
(2,0,0,0)a   &2 &6 &32 \\
(2,0,0,0)b   &2 &3 &8 \\
\hline
(1,1,1,1)a   &4 &4 & 96\\
(1,1,1,1)b   &4 &6 & 48\\
\hline
(1,1,1,0)   &3 &4 & 384\\
\hline
(1,1,0,0)a    &2 &7 & 24\\
(1,1,0,0)b    &2 &4 & 192\\
(1,1,0,0)c    &2 &5 & 192\\
(1,1,0,0)d    &2 &2 & 48\\
(1,1,0,0)e    &2 &3 &96 \\
\hline
(1,0,0,0)a    &1 &3 & 384\\
(1,0,0,0)b    &1 &1 & 96\\
(1,0,0,0)c    &1 &4 & 192\\
\hline
(0,0,0,0)a    &0 &0 &6 \\
(0,0,0,0)b    &0 &4 &96 \\
(0,0,0,0)c    &0 &4 & 48\\
(0,0,0,0)d    &0 &1 & 192\\
(0,0,0,0)e    &0 &2 & 96\\
(0,0,0,0)f    &0 &8 & 3\\
\end{array}
&\hspace{1in}\quad
\end{array}
\]
\vspace{-2mm}

\noindent
The total number of $\#E_{2,3}(g)\div 16$ in Table IV is $2520$.
Hence we can verify that $\#E_{2,3}=8!=40320$. 
We see that $c_{3}((2,0,0,0)b)=12$ and $c_{3}((1,1,0,0)e)=4$,
$c_{3}((0,0,0,0)b)=c_{3}((0,0,0,0)c)=0$,
$c_{4}((0,0,0,0)b)=14$ and $c_{4}((0,0,0,0)c)=12$.
In consequence, 
$SE_{2,3}=\coprod_{g\in {\goth g}_{2,3}}
SE_{2,3}(g)$.

\noindent
{\bf Acknowledgement:}
The authors would like to thank Yuzuru Maeda
for his discovery of two graphs, $(2,1,0)$ and $(2,2,1)$
in Table I of $\S$\ref{subsection:fourthtwo}.

\vspace{-5mm}
%
%

\end{document}